\documentclass{amsart}
  
  \usepackage{amsthm} 
  \usepackage{amsmath} 
  \usepackage{amssymb} 

\usepackage[tmargin=0.75in,bmargin=0.75in,lmargin=0.75in,rmargin=0.75in]{geometry}
\usepackage{float}
\usepackage{graphicx}
\usepackage{wrapfig}
\usepackage{lipsum}
\usepackage[percent]{overpic}
\usepackage{caption}
\usepackage{subcaption}
 \usepackage{graphicx}
  \usepackage{xspace}
  \usepackage[all]{xy}
  \usepackage{pinlabel}
  \usepackage{enumerate}
  \usepackage{hyperref}
\hypersetup{
    bookmarks=true,         
    unicode=false,          
    pdftoolbar=true,        
    pdfmenubar=true,        
    pdffitwindow=false,     
    pdfstartview={FitH},    
    pdftitle={My title},    
    pdfsubject={article},   
    pdfproducer={}, 
    pdfnewwindow=true,      
    colorlinks=true,       
    linkcolor=red,          
    citecolor=blue,        
    filecolor=magenta,      
    urlcolor=cyan           
}
\usepackage{tikz}
\usepackage{dsfont}
\usetikzlibrary{matrix}

  \usepackage{hyperref}

  \newcommand{\calM}{\mathcal{M}}

  \newcommand{\calQ}{\mathcal{Q}}

  \newcommand{\calT}{\mathcal{T}}

  \newcommand{\RR}{\mathbb{R}}

  \newcommand{\ZZ}{\mathbb{Z}}

  \theoremstyle{definition}

  \newtheorem*{claim*}{Claim}

  \newtheorem*{question*}{Question}
  \newtheorem*{answer*}{Answer}
  \newtheorem*{application*}{Application}
  \newtheorem*{sat*}{Theorem}

  \theoremstyle{remark}
  
  \newtheorem*{remark*}{Remark}
  
  \newcommand{\secref}[1]{Section~\ref{#1}}

  \newcommand{\figref}[1]{Figure~\ref{#1}}

  \newcommand{\eqnref}[1]{Equation~\eqref{#1}}

  \DeclareMathOperator{\Th}{Th}

 \DeclareMathOperator{\dd}{d}

  \newcommand{\ML}{\ensuremath{\mathcal{ML}}\xspace} 
   
  \newcommand{\Teich}{{Teichm\"uller }} 
  \newcommand{\param}{{\mathchoice{\mkern1mu\mbox{\raise2.2pt\hbox{$
  \centerdot$}}
  \mkern1mu}{\mkern1mu\mbox{\raise2.2pt\hbox{$\centerdot$}}\mkern1mu}{
  \mkern1.5mu\centerdot\mkern1.5mu}{\mkern1.5mu\centerdot\mkern1.5mu}}}

  \renewcommand{\setminus}{{\smallsetminus}}

\begin{document}
 \definecolor{adr}{cmyk}{0.99,0.,0.,0.1}
\newcommand{\ed}[1]{\textbf{\textcolor{adr}{#1}}}

\title[On normalizations  of  Thurston  measure  on  the  space  of  measured
laminations]
{On normalizations  of  Thurston  measure  on  the  space  of  measured
laminations}  

 \author   {Leonid Monin, Vanya Telpukhovskiy}
 \address{Department of Mathematics, University of Toronto, Toronto, ON, Canada}
 \email{lmonin@math.toronto.edu, ivantelp@math.toronto.edu}

\begin{abstract}  The space of measured laminations $\ML(\Sigma)$ associated to a topological surface $\Sigma$ of genus $g$ with $n$ punctures is an integral piecewise linear manifold of real dimension $6g-6+2n$. There is also a natural symplectic structure on $\ML(\Sigma)$ defined by Thurston. The integral and symplectic structures define a pair of measures on $\ML(\Sigma)$  which are known to be proportional. The projective class of these measures on $\ML(\Sigma)$ is called the Thurston measure. In this note we compute the ratio between two normailzations of the Thurston measure.
\end{abstract}
  
  \maketitle

\section{Introduction}

A real vector space $V$ has a unique up to scaling translation invariant measure $\mu$. 
Endowing $V$ with extra structures allows to provide natural normalizations of $\mu$. Here are some of them: 

\begin{itemize}
  \item  A Euclidean inner product $\langle\cdot,\cdot\rangle$ normalizes the volume by letting the volume of the parallelepiped spanned by the orthonormal basis be equal to 1;
  \item  A full rank lattice $\Lambda \subset V$, normalizes $\mu$ by letting $\mu(V/\Lambda)=1$;
  \item  If $\dim V=2n$, the standard symlpectic form $\omega$ defines a volume form $\frac{1}{n!}\omega^n$ and therefore normalizes $\mu$. 
\end{itemize}
  
In \cite{Masur} Masur showed that there exists a unique up to scaling mapping class group invariant measure on the space of measured laminations $\ML(\Sigma)$ in the Lebesgue class. As in the case of the vector space, there are different ways to normalize it, and we consider two
of them which are the most natural.

The first one comes from the symplectic form on $\ML(\Sigma)$ introduced by Thurston. In partucular, this normalization is natural due to the result of Bonahon and S\"ozen \cite{BS}, where they construct a family of symplectomorphisms between $\ML(\Sigma)$ and \Teich space $\calT(\Sigma)$ endowed with Weil-Petersson form.

The second normalization comes from the integral affine structure on $\ML(\Sigma)$.
It is natural for the counting problems on the space of quadratic differentials.
We call these normalizations (rigorously defined  in the \secref{back}) {\em symplectic} and {\em integral} Thurston measures on $\ML(\Sigma)$ and denote them by $ \mu_{\omega}, \, \mu_{\ZZ}$, respectively. In this paper we answer the question raised in \cite{Currents}: What is the ratio between measures $\mu_{\omega}$ and $\mu_{\ZZ}$? 

\begin{sat*}
The symplectic Thurston measure on $\ML(\Sigma)$ is $2^{|\chi(\Sigma)|-1}-$multiple of the integral Thurston measure:
$$ 
\frac{ \mu_{\omega} } { \mu_{\ZZ} } = 2^{|\chi(\Sigma)|-1} = 2^{2g+n-3}. 
$$
\end{sat*}
This result was also obtained in \cite{Arana-Herrera}. 
\subsection*{Related results.} 
The normalization of Thurston measure was a cause of the confusion in the literature before. One example is a theorem of Mirzakhani \cite{Maryam Ergodic}, that computes the Masur-Veech volume of the space of unit area quadratic differentials in terms of the integral of Mirzakhani function $B(X)$ over the WP volume form:
\begin{equation} \label{M formula}
Vol(\calQ^1\calM_g) = \int_{\calM_g} B(X) \, dX.
\end{equation}
Although this formula is essentially correct, if one chooses the standard conventions for Thurston measure and WP volume, the formula requires additional constant. Our theorem computes essential part of this constant. For more details, see \cite{DGZZ19}, where a version of \eqnref{M formula} with the correct normalization was also obtained.

\subsection*{Acknowledgements.} We would like to thank Kasra Rafi for introducing the counting problems on \Teich space to us. We also thank Kasra Rafi and Anton Zorich for their interest and support.

\section{Background}\label{back}

Fix a topological surface $\Sigma$ of genus $g$ with $n$ punctures with negative Euler characteristic. Let $\calT(\Sigma)$ be the \Teich space of the surface $\Sigma$, i.e. the space of isomorphism classes of marked hyperbolic metrics on $\Sigma$. Let $\ML(\Sigma)$ be the space of measured laminations on $\Sigma$. For the introduction to these geometric structures we refer the reader to \cite{Bruno}.

\subsection{Train track coordinates}
A \emph{train track} on $\Sigma$ is an embedded 1-complex $\tau$ such that 
\begin{itemize}
    \item each edge (branch) of $\tau$ is a smooth path with well-defined tangent vectors at the end points. That is, all edges at a given vertex (switch)
are tangent to each other. 
    \item For each component $R$ of $\Sigma \setminus \tau$, the double of $\partial R$ along the interior of the edges of $R$ has negative Euler characteristic.
\end{itemize}

For a train track $\tau$ we will denote the set of its branches by $b(\tau)$ and the set of its switches by $sw(\tau)$. A (measured) lamination $\lambda$ is \emph{carried} by $\tau$ if there is a differentiable map $f: \Sigma \to \Sigma$ homotopic to identity taking $\lambda$ to $\tau$, such that the restriction of $df$ to tangent lines of $\lambda$ is non-singular. 

A measured lamination $\lambda$ carried by the train track $\tau$ assigns the weight to each branch of $\tau$. The weights coming from a measured lamination are non-negative and satisfy the {\em switch conditions}: for every switch the sums of the weights of the incoming and outgoing branches are equal. Moreover, each collection of positive weights satisfying switch conditions correspond to some measured lamination carried by $\tau$. We will denote the subspace of the vector space $\RR^{b(\tau)}$ satisfying switch conditions by $W(\tau)$. Therefore, the set $V(\tau)$ of measured laminations that are carried by $\tau$ is identified with a polyhedral cone in $W(\tau)$, so that $V(\tau)\cong W(\tau)\cap\RR_+^{b(\tau)}$.

A train track $\tau$ is called \emph{recurrent} if each branch $b$ of $\tau$ is traversed by some curve $\alpha$, carried by $\tau$. A train track $\tau$ is called \emph{transversely recurrent} if for each branch $b$ of $\tau$ there is a dual curve $\gamma$ intersecting $b$. A train track $\tau$ is called \emph{birecurrent} if it is both recurrent and transversely recurrent. A train track is called {\em maximal} if its complement $\Sigma \setminus \tau$ is a union of trigons and once punctured monogons. The train track is called \emph{generic} if it has only trivalent switches. Since every lamination is carried by some generic maximal birecurrent train track, the the collection of sets $V(\tau)$ for such train tracks gives rise to an atlas of charts on $\ML(\Sigma)$.

\subsection{Integral lattice and the measure $\mu_\ZZ$}

We will denote by $\Lambda(\tau)$ the lattice of integral points in $W(\tau)$, that is $\Lambda(\tau)=W(\tau)\cap\ZZ^{b(\tau)}$. Geometrically  $\Lambda(\tau)$ corresponds to integrally weighted multicurves carried by $\tau$. For a pair of train tracks $\tau, \tau'$ the transition map between the cones $V(\tau)$ and $V(\tau')$ is an integral linear transformation with respect to the lattices $\Lambda(\tau)$, $\Lambda(\tau')$. Thus the train track charts define an integral piecewise linear structure on $\ML(\Sigma)$.

A collection of lattices $\Lambda(\tau)$ for all birecurrent train tracks define normalizations of a volume form on $W(\tau)$, and therefore on $V(\tau)$. These normalizations are consistent since the transition maps in train track coordinates preserve the lattices. A global volume form defined in such a way provides the measure  $\mu_\ZZ$.

{\makeatletter
\let\par\@@par
\par\parshape0
\everypar{}
\begin{wrapfigure}{r}{0.28\textwidth}
     \centering
     \begin{overpic}[scale=0.2]{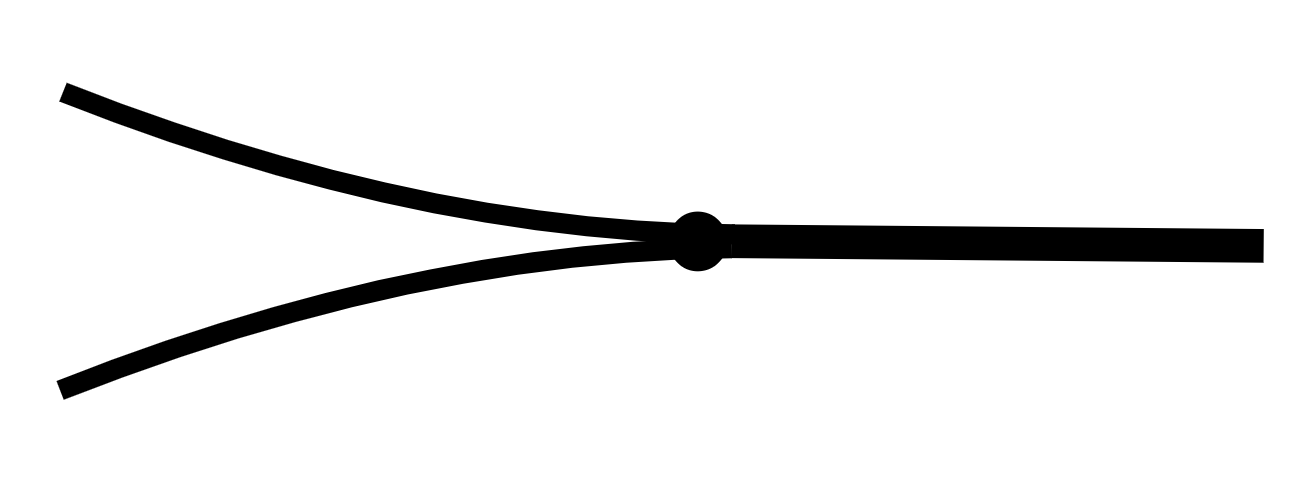}
 \put (20,29) {$s_1$}
 \put (20,8) {$s_2$}
 \put (52,11) {$s$}
\end{overpic}
\caption{}
\label{switch}
\end{wrapfigure}
\subsection{Symplectic form and the measure $\mu_\omega$}
For a point $\lambda\in V(\tau)$, the tangent space $T_\lambda\ML$ is naturally identified with the vector space $W(\tau)$. The symplectic structure on $\ML(\Sigma)$ is defined by the collection of skew symmetric bilinear forms $\omega(\cdot,\cdot)_\tau$ for all generic maximal birecurrent train tracks $\tau$, that are consistent with respect to the transition maps. Such a system was defined by Thurston as a restriction of a skew symmetric form $\widetilde \omega$ on $\RR^{b(\tau)}$, where:\par}
$$
\widetilde \omega(x,y)_\tau =\frac{1}{2}\sum_{s\in sw(\tau)} s_1(x)s_2(y)-s_2(x)s_1(y),
$$ 
where $s_1, s_2$ are two incoming edges of a switch $s$ as on \figref{switch}. In other terms, the Thurston symplectic form in the chart $V(\tau)$ is given by:
$$
\widetilde\omega=\frac{1}{2}\sum_{s\in sw(\tau)} \dd s_1\wedge \dd s_2.
$$
It is a theorem of Thurston that the restrictions of $\widetilde \omega$ to $W(\tau)$ is a collection of nondegenerate bilinear pairings which define a mapping class group invariant symplectic form $\omega_{\Th}$ on $\ML(\Sigma)$. The natural volume form defined by  $\frac{1}{(3g-3+n)!}\omega_{\Th}^{3g-3+n}$   gives rise to a measure on $\ML(\Sigma)$ that we denote by $\mu_\omega$.

\subsection{Euclidean normalization of Thurston measure is not well defined} 
Another potential normalization of Thurston measure comes from the standard Euclidean structure on the space $\RR^{b(\tau)}$.
For a given train track $\tau$, the standard Euclidean structure on $\RR^{b(\tau)}$ defines a volume form on the space $W(\tau)$ of weights satisfying switch conditions, and therefore on $V(\tau)$.
It is proportional to Thurston measure locally in every chart $V(\tau)$. Nevertheless, the Euclidean normalization is not globally well-defined. 

\begin{figure}[H]
\centering
\begin{subfigure}{.5\textwidth}
  \centering
  \includegraphics[scale=0.3]{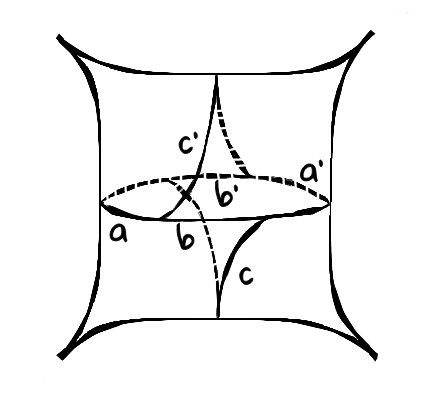}
  \caption{Train track $\tau$ on $\Sigma_{0,4}$}
\end{subfigure}%
\begin{subfigure}{.5\textwidth}
  \centering
 \includegraphics[scale=0.3]{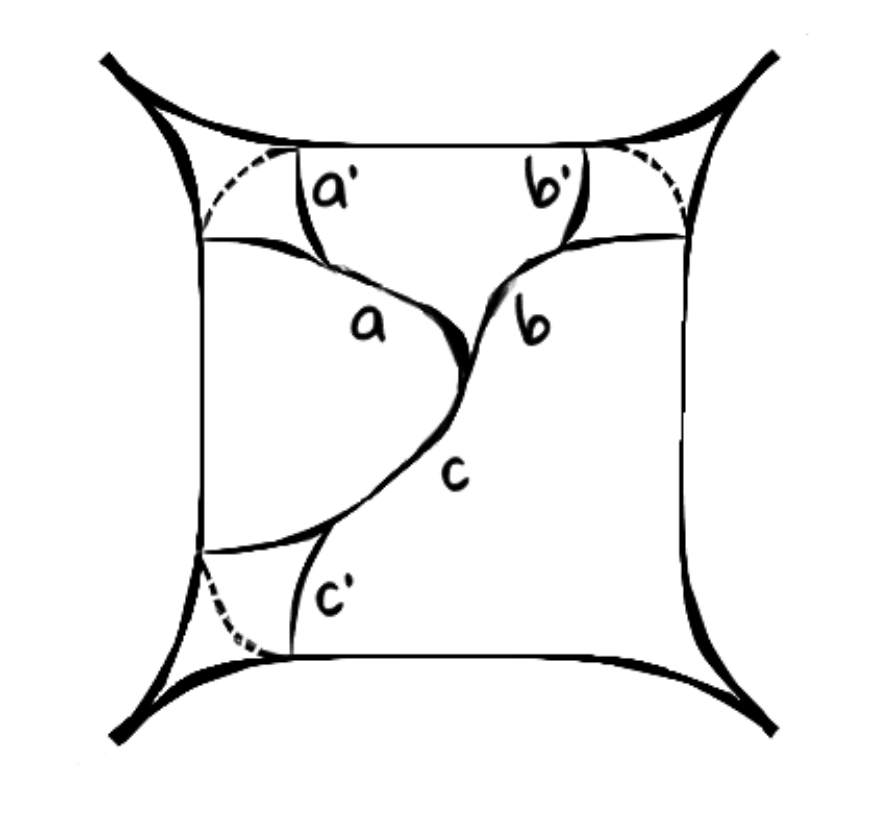}
    \caption{Train track $\tau'$ on $\Sigma_{0,4}$}
\end{subfigure}
\caption{Two train tracks corresponding to different Euclidean normalizations of $\mu_{\Th}$.}\label{Euclidlox}
\end{figure}

To show this we consider two maximal birecurrent train tracks $\tau$, $\tau'$ on the four-times punctured sphere as on Figure~\ref{Euclidlox} and compare corresponding Euclidean normalizations with the integral normalization. 
\begin{enumerate}
\item  For the train track $\tau$, we have $\RR^{b(\tau)}=\langle a,b,c,a',b',c'\rangle$ and the switch conditions are given by the following equations: 
$$
a = b+c' = b'+c,\quad a'= b+c = b'+c'.
$$
For the basis of the lattice $\Lambda(\tau)$ given by the vectors $v_1=(1,1,0,1,1,0)$, $v_2=(1,0,1,1,0,1),$ the Euclidean area of its fundamental parallelepiped equals to $2\sqrt{3}$.

\item  For the train track $\tau'$, we have $\RR^{b(\tau')}=\langle a,b,c,a',b',c'\rangle$ and the switch conditions are given by the equations: 
$$
a =2a',\quad b=2b',\quad c=2c', \quad a+b=c.
$$
For the basis of the lattice $\Lambda(\tau')$ given by the vectors $v_1=(2,0,2,1,0,1)$, $v_2=(0,2,2,0,1,1),$ the Euclidean area of its fundamental parallelepiped equals to $5\sqrt{3}$.
\end{enumerate}

\noindent
Since the Euclidean areas of the fundamental parallelograms of the lattices $\Lambda(\tau), \Lambda(\tau')$ are not equal, the global Euclidean normalization is not well-defined.


\section{Proof of the main theorem} To find the ratio between measures $\mu_\omega$ and $\mu_\ZZ$ it is enough to do the computation locally. More precisely, it is enough to compute the volume of the torus $W(\tau)/\Lambda(\tau)$ with respect to the volume form $\frac{1}{(3g-3+n)!}\omega_{\Th}^{3g-3+n}$ for some maximal birecurrent train track $\tau$. We define such a train track $\tau_{g,n}$ for each pair $\{g, n\}$, so that there is a natural embedding of the space $W(\tau_{g,n})$ into both spaces $W(\tau_{g+1,n})$ and $W(\tau_{g,n+1})$. We denote the symplectic form $\omega_{\Th}$ on the space $W(\tau_{g,n})$ by $\omega_{g,n}$. We show that under the above embedding, the expressions $\omega_{g+1,n}-\omega_{g,n}$ and $\omega_{g,n+1}-\omega_{g,n}$ have  standard forms (see the \secref{expressions for forms}). We use this to prove our main theorem by induction.

\subsection{Induction scheme} 
\label{subsec: ind sch}
We prove the main result using the induction with respect to two parameters: the genus and the number of punctures on the surface. For closed and once punctured surfaces we run induction by adding a genus with base cases $\Sigma_{2,0}$, $\Sigma_{2,1}$ (see the first two columns in the \figref{induction}). For the remaining surfaces we run induction by adding a puncture using surfaces $\Sigma_{0,5}$, $\Sigma_{1,2}$, and $\Sigma_{g,1}$ for $g\geqslant 2$ as base cases. We leave the remaining three cases of $\Sigma_{0,3}$, $\Sigma_{0,4}$, and $\Sigma_{1,1}$ for the reader. Here is the full scheme of the induction:
\begin{figure}[H]
\begin{center}
\begin{tikzpicture}
  \matrix (m) [matrix of math nodes,
    nodes in empty cells,nodes={minimum width=4ex,
    minimum height=4ex,outer sep=-4pt},
    column sep=0.8ex,row sep=0.8ex]{
          g     &      &     &     &    &\vdots&\vdots&     &     & \\
                &      &     &     &     &      &     &     &     & \\
          2     &   B  &  B  &     &     &      &     &     &     & \\
          1     &   *  &  B  &  B  &     &      &     &     &     & \\
          0     &   *  &  *  &  *  &  B  &   B  &  B  &     &     & \\
    \quad\strut &   0  &  1  &  2  &  3  &   4  &  5  &     &  n  & \strut \\};
  \draw[-stealth] (m-2-2.south) -- (m-1-2.north);
  \draw[-stealth] (m-2-3.south) -- (m-1-3.north);
  \draw[-stealth] (m-3-4.center) -- (m-3-10); 
  \draw[-stealth] (m-2-4.center) -- (m-2-10);
  \draw[-stealth] (m-4-5.center) -- (m-4-10);
  \draw[-stealth] (m-5-8.center) -- (m-5-10);
  \draw[thick] (m-6-1.east) -- (m-1-1.east) ;
  \draw[thick] (m-6-1.north) -- (m-6-9.north) ;
\end{tikzpicture}
\end{center}
\caption{The induction scheme.}\label{induction}
\end{figure}
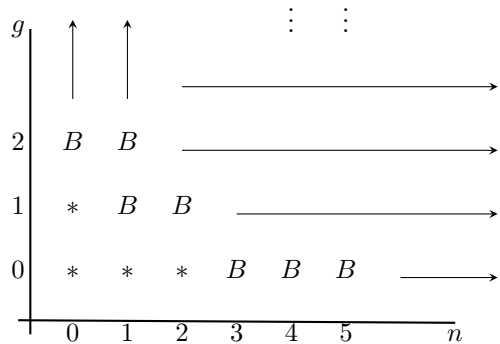
By * we denote the surfaces with non-negative Euler characteristic that aren't considered in the paper and by $B$ we denote the base cases of the induction (see the \secref{subsection: base cases} for the computation). Directions of the arrows represent the directions of the induction: adding genus or the puncture. 

\subsection{Base cases} 
\label{subsection: base cases}
In this subsection we consider the base cases of the inductions described in the \secref{subsec: ind sch}.

\subsection*{Five-times punctured sphere}
In the case of the five-times punctured sphere, we work in the chart defined by the train track $\tau_{0,5}$ as on \figref{fig:05 and 12}. The weight on each branch of $\tau_{0,5}$ can be expressed as a linear combination of the weights $z_1, \dotsc, z_4$ using the switch conditions. This defines an isomorphism between spaces $\RR(z_1, \dotsc, z_4)$ and $W(\tau_{0,5})$. It is also easy to see that this  is an isomorphism on the level of the integer lattices $\ZZ(z_1, \dotsc, z_4)$ and $\Lambda(\tau_{0,5})$. In the coordinates $z_1, \dotsc, z_4$, the expression for Thurston symplectic form in this chart is given by
$$\omega_{0,5} = \frac{1}{2} \left( 4 \dd z_1\wedge\dd z_2 + 4 \dd (z_1 + z_2)\wedge\dd z_3 + 4 \dd (z_1+z_2+z_3)\wedge\dd z_4\right). $$
We compute the volume form, which is given by $\frac{\omega_{0,5}^2}{2!} = 4 \dd z_1\wedge\dd z_2\wedge\dd z_3\wedge\dd z_4.$
The integer lattice $\Lambda(\tau_{0,5})$
coincides with the standard integer lattice in the coordinates $z_1,\ldots, z_4$.
Therefore the ratio between measures $\mu_{\omega}$ and $\mu_{\ZZ}$ in this case equals to
$4 = 2^{|\chi(\Sigma_{0,5})|-1} $.

\subsection*{Twice punctured torus}
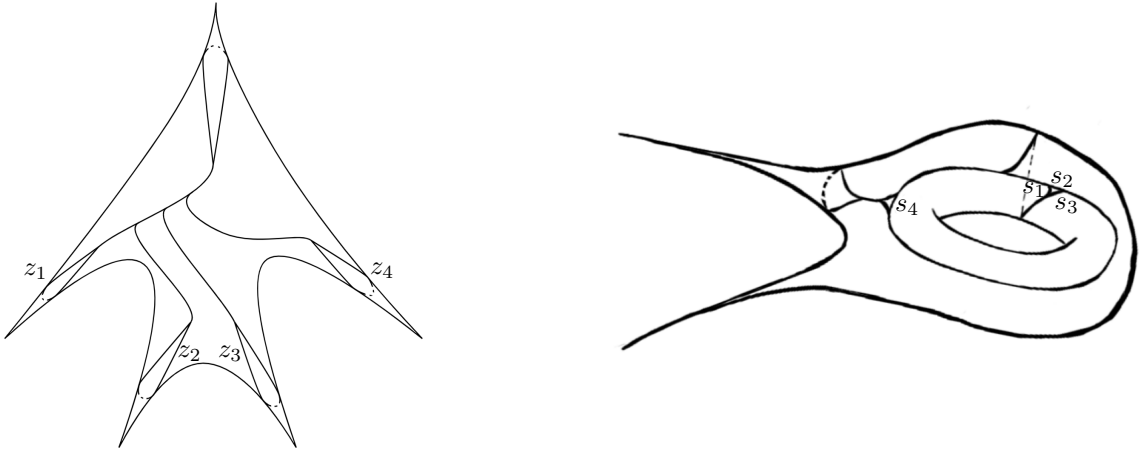
\begin{figure}[H]
\centering
\begin{subfigure}{.5\textwidth}
  \centering
  \begin{overpic}[scale=0.32]{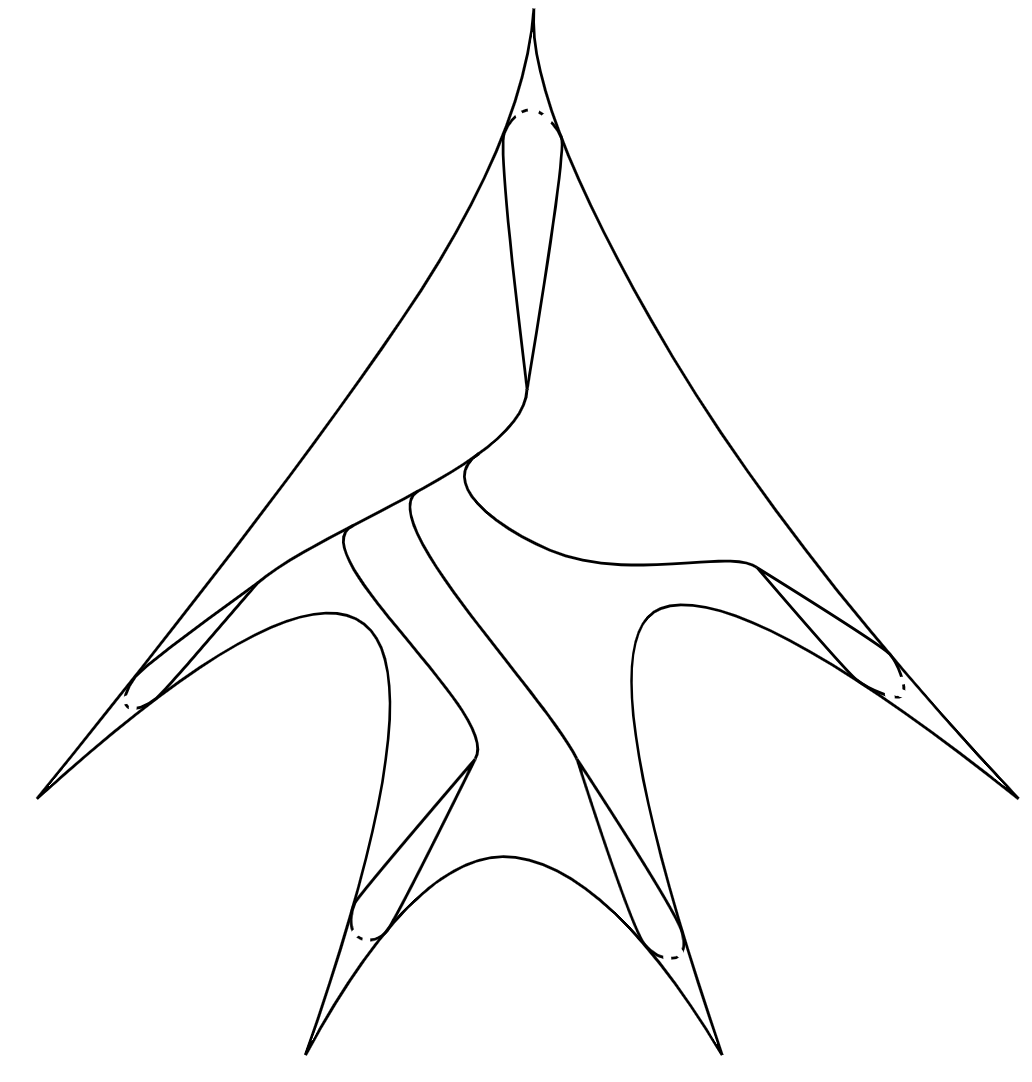}
 \put (7.5,40.) {$z_1$}
 \put (41,22.) {$z_2$}
 \put (49.9,22) {$z_3$}
 \put (83,40.1) {$z_4$}
\end{overpic}
  \label{fig:sub1}
\end{subfigure}%
\begin{subfigure}{.5\textwidth}
  \centering
  \begin{overpic}[scale=0.25]{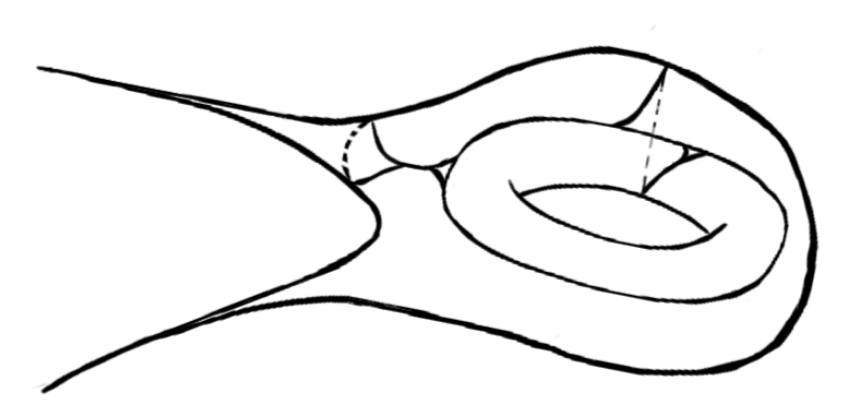}
 \put (76.0,30.3) {$s_1$}
 \put (81.,32.2) {$s_2$}
 \put (81.8,27.5) {$s_3$}
 \put (53.4,27.2) {$s_4$}
\end{overpic}
  \label{fig:sub2}
\end{subfigure}
\caption{Train tracks $\tau_{0,5}$ and $\tau_{1,2}$.}\label{m0npic}
\label{fig:05 and 12}
\end{figure}

In the case of the twice punctured torus, we work in the chart defined by the train track $\tau_{1,2}$ as on \figref{fig:05 and 12}. The weight of every branch of $\tau_{1,2}$ can be expressed through the weights  $s_1,\dotsc,s_4$. This, as in the case of the five-times punctured sphere, defines the coordinates on the space $W(\tau_{1,2})$. In these coordinates, Thurston symplectic form is given by
\begin{equation*}
\label{equation: 12}
\begin{split}
\omega_{1,2} = \frac{1}{2} ( \dd s_1\wedge\dd s_2 + \dd s_3\wedge\dd s_1 + \dd s_2\wedge\dd s_3 + \dd (s_1+s_2)\wedge\dd (s_1+s_3) + \dd s_4\wedge\dd (s_2+s_3) \\+ \dd (2s_1+s_2+s_3)\wedge\dd s_4 + \dd (s_2+s_3-s_4) \wedge\dd (2s_1+s_2+s_3-s_4)).
\end{split}
\end{equation*}
The volume form then equals to $\frac{\omega_{1,2}^2}{2!} = 2 \dd s_1\wedge\dd s_2\wedge\dd s_3\wedge\dd s_4.$
The lattice  $\Lambda(\tau_{1,2})$ is the standard integer lattice in the coordinates $s_1,\ldots,s_4$, so the ratio between measures $\mu_{\omega}$ and $\mu_{\ZZ}$ equals to $2 = 2^{|\chi(\Sigma_{1,2})|-1}.$

\subsection*{Genus two surface}
In the case of the genus two surface, we work in the chart given by the train track $\tau_{2,0}$ as on \figref{fig:21 and 20} (right), and the coordinates on the space $W(\tau_{2,0})$ are given by the weights on the branches $x_1,\ldots,x_6$. In these coordinates, the
Thurston symplectic form is given by $$\omega_{2,0} = \frac{1}{2}\left(2\dd x_1\wedge \dd x_2 +  2\dd (2x_1 + 2x_2 - 2x_3)\wedge \dd x_4  +  2\dd x_5 \wedge \dd x_6\right).$$ The volume form equals to
$
\frac{\omega_{2,0}^3}{3!} = \pm 2 \dd x_1\wedge\dd x_2\wedge\dd x_3\wedge\dd x_4
\wedge\dd x_5\wedge\dd x_6.
$
It easy to see that if the weights on the branches $x_1,\ldots, x_6$ are integer then all other weights are intefer too. The lattice $\Lambda(\tau_{2,0})$ is the standard lattice in the coordinates $x_1,\ldots, x_6$, therefore the ratio between symplectic and integral Thurston measures equals to $2 = 2^{|\chi(\Sigma_{2,0})|-1}$.

\subsection*{Genus two surface with a puncture}
\label{subsection: g21}
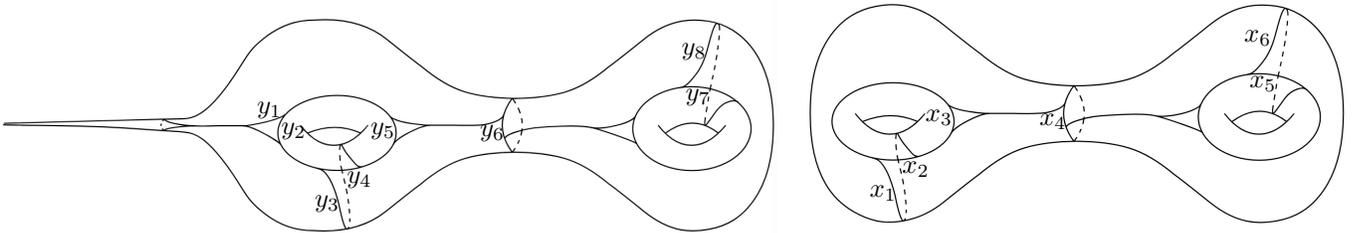
\begin{figure}[H]
\centering
\begin{subfigure}{.6\textwidth}
  \centering
  \begin{overpic}[scale=0.4]{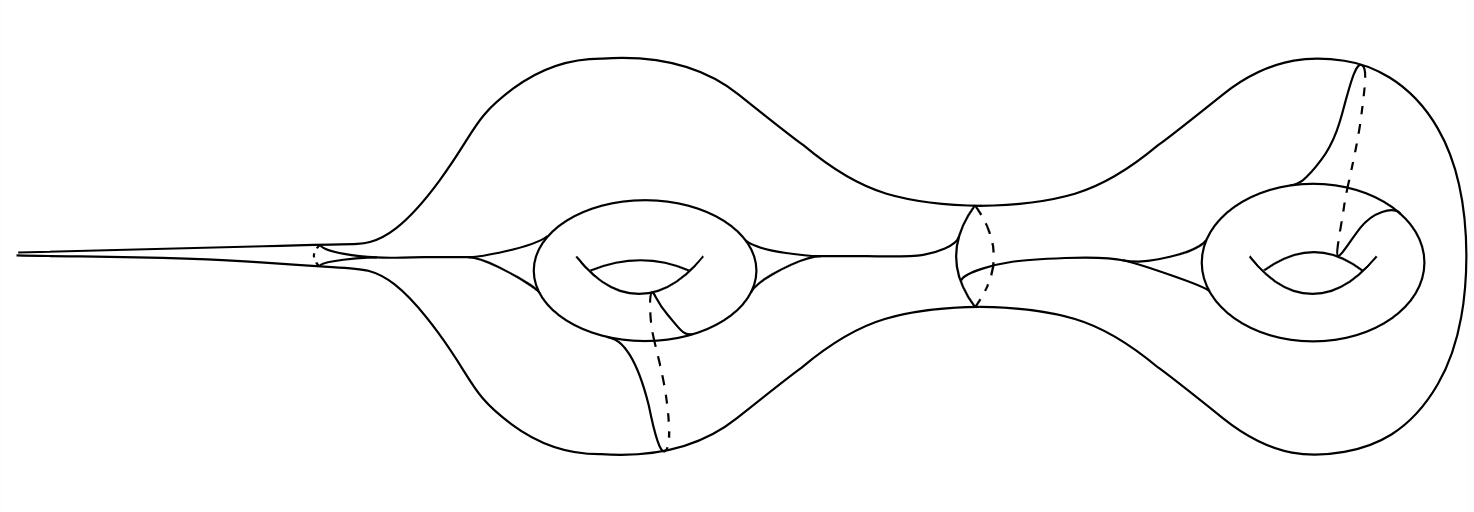}
 \put (33.5,18.9) {$y_1$}
  \put (36.6,16.2) {$y_2$}
 \put (40.9,7) {$y_3$}
 \put (45,10) {$y_4$}
 \put (48,16.2) {$y_5$}
 \put (62.0,16) {$y_6$}
 \put (88.3,20.6) {$y_7$}
 \put (87.8,26.4) {$y_8$}
\end{overpic}
  \label{fig:sub1}
\end{subfigure}%
\begin{subfigure}{.4\textwidth}
  \centering
  \begin{overpic}[scale=0.3]{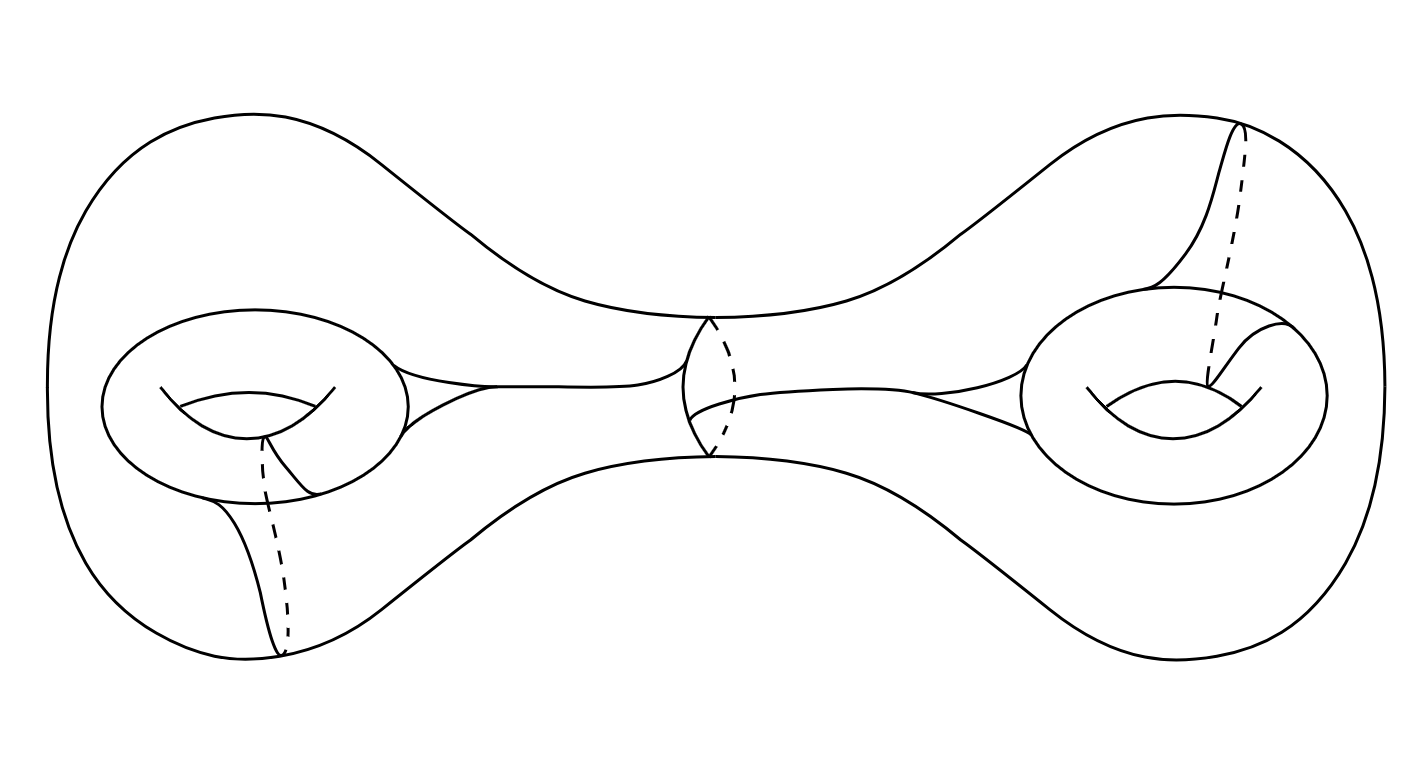}
 \put (13.8,11.5) {$x_1$}
 \put (19.6,15.7) {$x_2$}
 \put (23.7,25) {$x_3$}
 \put (43.9,24.2) {$x_4$}
 \put (81.1,31.0) {$x_5$}
 \put (80.2,39) {$x_6$}
\end{overpic}
  \label{fig:sub2}
\end{subfigure}
\caption{Train tracks $\tau_{2,1}$ and $\tau_{2,0}$.}
\label{fig:21 and 20}
\end{figure}

In the case of the genus two surface with a puncture, we work in the chart  given by the train track $\tau_{2,1}$ as on \figref{fig:21 and 20} (left), and the coordinates on the space $W(\tau_{2,0})$ are given by the weights on the branches $y_1,\ldots,y_8$. In these coordinates, the
Thurston symplectic form is given by
\begin{equation*}
\begin{split}
\omega_{2,1} = \frac{1}{2}\Big(\dd y_1 \wedge\dd y_2 + \dd y_2\wedge\dd (y_3+  y_4)+\dd(y_3+y_4-y_2)\wedge \dd y_1 + 2\dd y_3\wedge\dd y_4 +\dd(y_3+y_4)\wedge\dd y_5 +\\
\dd y_5\wedge \dd (y_1+y_2)+ \dd(y_1+y_2-y_5)\wedge \dd(y_3+y_4-y_5)+2\dd(y_1+y_2+y_3+y_4-2y_5)\wedge\dd y_6+ 2\dd y_7\wedge\dd y_8\Big).
\end{split}
\end{equation*}
The volume form equals to
$
\frac{\omega_{2,1}^4}{4!} = \pm 2 \dd y_1\wedge\dd y_2\wedge\dd y_3\wedge\dd y_4
\wedge\dd y_5\wedge\dd y_6 \wedge\dd y_7 \wedge\dd y_8.
$
It easy to check that the lattice $\Lambda(\tau_{2,1})$ is given by the integer weights on the branches $y_1,\ldots, y_8$, such that the sum $y_1-y_2+y_3+y_4$ is even. In other words, the lattice $\Lambda(\tau_{2,1})$ is the kernel of the homomorphism $\phi:\ZZ^8\to \ZZ/2\ZZ$ given by 
$$
\phi(y_1,\ldots, y_8)=y_1-y_2+y_3+y_4 \,\, (mod \,\, 2),
$$
where $\ZZ^8$ is the standard lattice in the coordinates $y_1,\ldots, y_8$. Therefore, the lattice $\Lambda(\tau_{2,1})$ is the sublattice in $\ZZ^8$ of index $2$, and the ratio between symplectic and integral Thurston measures is equal to $4 = 2^{|\chi(\Sigma_{2,1})|-1}$.

\subsection{Train tracks $\tau_{g,n}$ and the expressions for $\omega_{g,n}$}
\label{subsec: descriptions}
\subsection*{Train tracks $\tau_{g,n}$.}
\begin{figure}
  \begin{overpic}[scale=0.215]{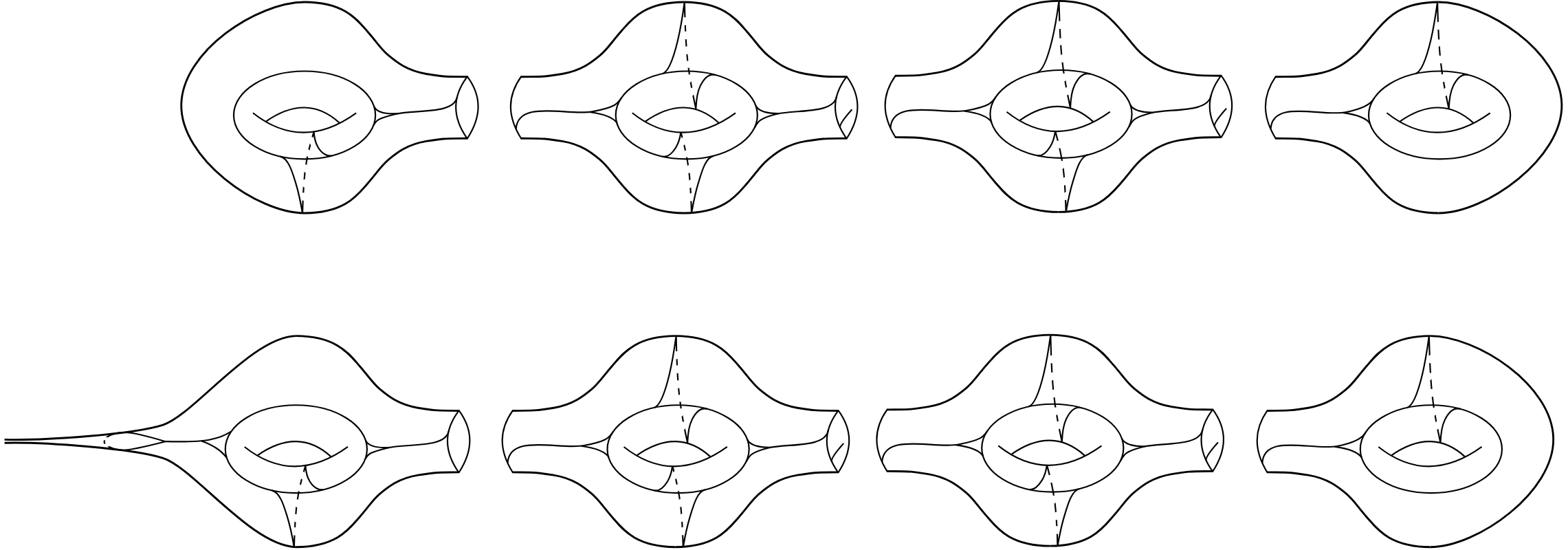}
\put (17,18.3) {$X_{left}$}
\put (17,15.7) {$Y_{left}$}
\put (43.5,17) {$S_1$}
\put (55,17.) {$\cdots$}
\put (67.5,17.) {$S_{g-2}$}
\put (90.5,18.3) {$X_{right}$}
\put (90.5,15.7) {$Y_{right}$}
\end{overpic}
\caption{Train tracks $\tau_g$ and $\tau_{g,1}$.}\label{gg1}
\end{figure}
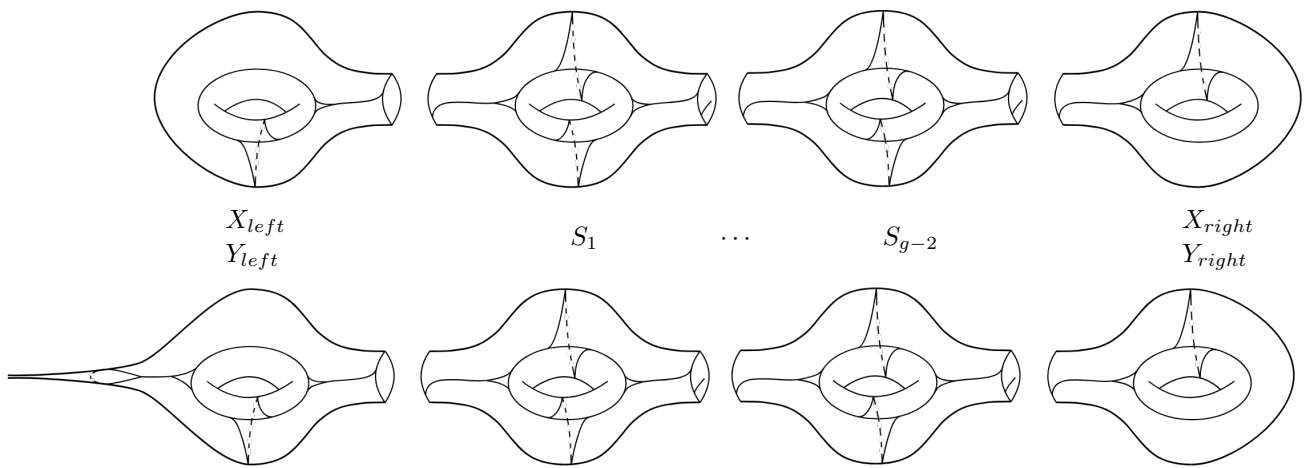

First, we define the train tracks $\tau_{g,0}$ on closed surfaces of genus at least 2. We decompose the surface into $g$ subsurfaces, such that $g-2$ of them are homeomorphic to a torus with two boundary components denoted by $S_1, \dotsc, S_{g-2}$ and the other two $X_{left},X_{right}$ are homeomorphic to a torus with one boundary component, as on Figure~\ref{gg1}. The restrictions of the train track $\tau_{g,0}$ to the subsurfaces $S_1,\dotsc,S_{g-2}$ are the same and as on Figure~\ref{fig: triangle and square} (right). The restriction to the remaining two subsurfaces $X_{left}, X_{right}$ can be glued to a train track $\tau_{2,0}$ as on Figure~\ref{fig:21 and 20}.

For the train tracks $\tau_{g,1}$ on once punctured surfaces of genus at least 2, we proceed similarly. We decompose the surface into $g$ subsurfaces, such that $g-2$ of them are homeomorphic to a torus with two boundary components denoted by $S_1, \dotsc, S_{g-2}$, the subsurface $Y_{left}$ homeomorphic to a punctured torus with one boundary component, and the subsurface $Y_{right}$  homeomorphic to a torus with one boundary component, as on Figure~\ref{gg1}. The restrictions of the train track $\tau_{g,1}$ to $g-2$ tori with two boundary components $S_1, \dotsc,S_{g-2}$ are as before, see Figure~\ref{gg1}. The restriction to the remaining two subsurfaces $Y_{left},Y_{right}$ can be glued to a train track $\tau_{2,1}$ as on Figure~\ref{fig:21 and 20}.

Denote by $\Sigma_{g,n}^k$ the surface of genus $g$ with $n$ punctures and $k$ boundary components. To define the general train track $\tau_{g,n}$ we consider three cases. In each of them we cut the surface into two subsurfaces by a simple closed curves, such that one subsurface is $\Sigma_{0,n}^1$ and the restriction of the train track $\tau_{g,n}$ to this subsurface is as on Figure~\ref{fig:cactus}. 

In the case of the sphere with $n$ punctures ($n\geqslant 5$), we decompose the sphere into subsurface $\Sigma_{0,n-4}^1$ and the complementary subsurface $C$ homeomorphic to $\Sigma_{0,4}^1$. We complete the train track $\tau_{0,n}$ in such a way that the red subtrack on Figure~\ref{fig:cactus} together with the restriction of $\tau_{0,n}$ to the subsurface $C$ form the train track $\tau_{0,5}$ as on \figref{fig:05 and 12}.

In the case of the torus with $n$ punctures ($n \geqslant 2$), we decompose the surface into subsurface $\Sigma_{0,n-1}^1$ and the complementary subsurface $C$ homeomorphic to $\Sigma_{1,1}^1$. We complete the train track $\tau_{1,n}$ in such a way that the red subtrack on Figure~\ref{fig:cactus} together with the restriction of $\tau_{1,n}$ to the subsurface $C$ form the train track $\tau_{1,1}$.

In the case of the genus $g$ surface with $n$  punctures ($g \geqslant 2, n \geqslant 1$), we decompose the surface into subsurface $\Sigma_{0,n}^1$ and the complementary subsurface $C$ homeomorphic to $\Sigma_{g,0}^1$. We complete the train track $\tau_{g,n}$ in such a way that the red subtrack on Figure~\ref{fig:cactus} together with the restriction of $\tau_{g,n}$ to the subsurface $C$ form the train track $\tau_{g,1}$.

\begin{center}
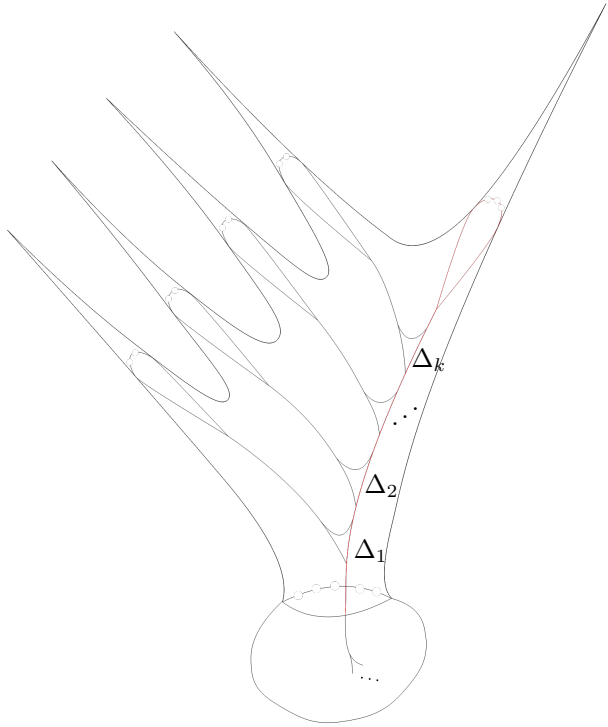
\begin{figure} 
\begin{overpic}[scale=0.15]{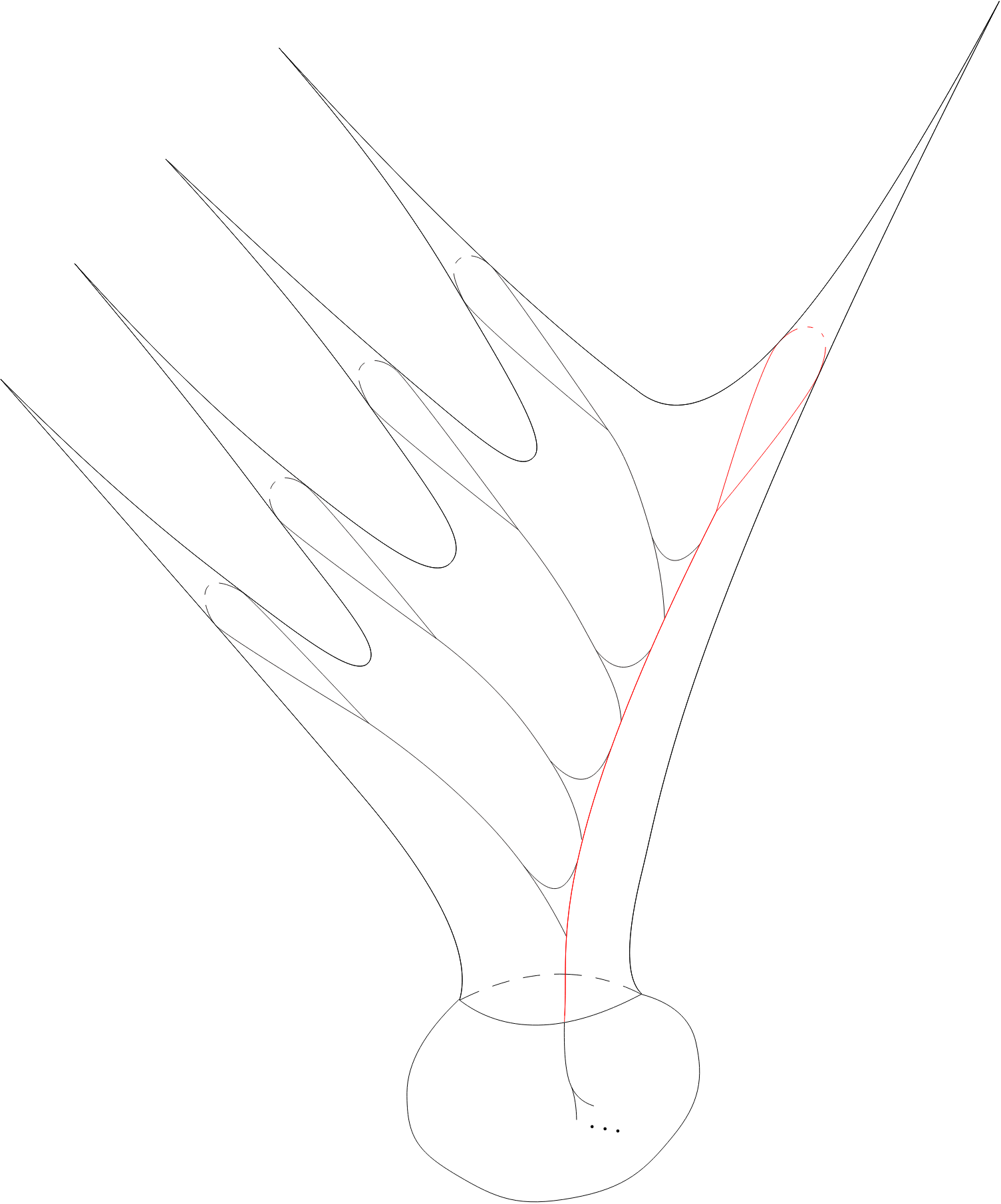}
\put (48,23) {$\Delta_1$}
\put (49.5,32) {$\Delta_2$}
\put (53,41) {\reflectbox{$\ddots$}}
\put (56,49.5) {$\Delta_k$}
\end{overpic}
\caption{On top: the subsurface $\Sigma_{0,n}^1$ with the restriction of the train track $\tau_{g,n}$ to it. \\\hspace{\textwidth}
On the bottom: the complementary subsurface $C$.}
\label{fig:cactus}
\end{figure}
\end{center}

\subsection*{Expressions for $\omega_{g,n}$.}
\label{expressions for forms}
Having defined the train tracks $\tau_{g,n}$, we proceed by expressing inductively the symplectic forms $\omega_{g,n}$  following the induction scheme in the \secref{subsec: ind sch}. As it was mentioned in the \secref{subsec: ind sch}, for closed and once punctured surfaces of genus at least $2$, the base cases of the induction are the surfaces $\Sigma_{2,0}$ and $\Sigma_{2,1}$. 

\begin{figure}[H]
\centering
\begin{subfigure}{.5\textwidth}
  \centering
  \begin{overpic}[scale=0.25]{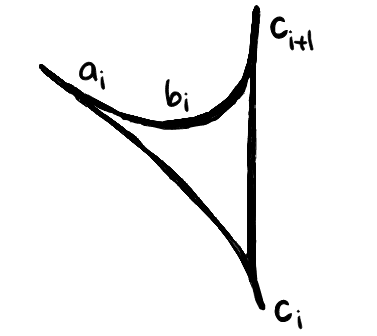}
 \put (14,68) {$2$}
\end{overpic}
  \label{fig:sub1}
\end{subfigure}%
\begin{subfigure}{.5\textwidth}
  \centering
  \begin{overpic}[scale=0.22]{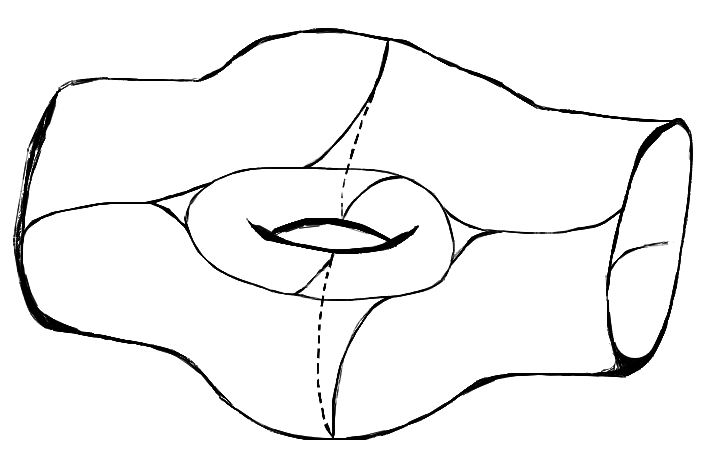}
 \put (19,41) {$A_i$}
 \put (27.9,34.9) {$B_i$}
 \put (44.9,55) {$C_i$}
 \put (49.,10.1) {$D_i$}
  \put (57,31) {$E_i$}
 \put (78,25.9) {$F_i$}
 \put (5,30) {$G_i$}
  \put (87.9,26.5) {$G_{i+1}$}
\end{overpic}
  \label{fig:sub2}
\end{subfigure}
\caption{The trigon $\Delta_i$ and the torus with two boundary components $\Box_i$.}\label{triangle and square}
\label{fig: triangle and square}
\end{figure}
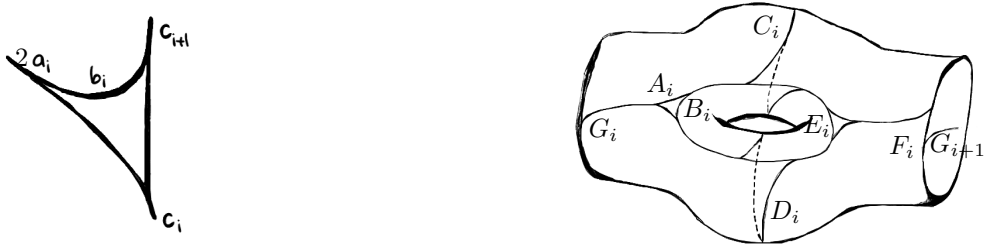

To define the coordinates on the space $W(\tau_{g,0})$ we use the branches $A_i,\ldots, F_i$ on subsurfaces $S_i$ as on \figref{triangle and square} (right) and branches $x_1,\ldots,x_6$ on the remaining two subsurfaces $X_{left}, X_{right}$ as on Figure~\ref{fig:21 and 20} (right).

Similarly, to define the coordinates on the space $W(\tau_{g,1})$ we also use the branches $A_i,\ldots, F_i$ on subsurfaces $S_i$ and the branches $y_1,\ldots,y_8$ in the remanining two subsurfaces $Y_{left}, Y_{right}$ as on Figure~\ref{fig:21 and 20} (left).

In these cases, the symplectic forms can be written as follows: 
\begin{equation} 
\label{equation: g0g1}
\omega_{g,0} = \omega_{2,0} +\sum_{i=1}^{g-2} \Box_i, \,\,\,\,\,\,\,\,\,\,
\omega_{g,1}= \omega_{2,1} +\sum_{i=1}^{g-2} \Box_i,
\end{equation}
where $\Box_i$ is the contribution of the switches on the subsurface $S_i$ and the expressions for the symplectic forms $\omega_{2,0}$ and $\omega_{2,1}$ are as in the \secref{subsection: base cases}. We note that this choice of coordinates allows us naturally embed spaces $W(\tau_{g,0}), W(\tau_{g,1})$ into the spaces $W(\tau_{g+1,0}), W(\tau_{g+1,1})$ and under this embedding 
\begin{equation}
\label{eq: recursive puncture}
\omega_{g+1,0}= \omega_{g,0} +\Box_{g-1}, \quad \omega_{g+1,1}= \omega_{g,1} +\Box_{g-1}.
\end{equation}
To give an explicit expression for $\Box_i$ we introduce auxiliary variables given by weights of branches $G_i$ on the subsurfaces $S_i$ as on \figref{triangle and square} (right). The weights of branches $G_i$ can be expressed  through the other variables inductively in the following way:
$$
G_{i+1} = G_i + 2B_i -2E_i,
$$
with $G_1=x_1+x_2-2x_3$ for closed surfaces and $G_1=y_1+y_2+y_3+y_4-2y_5$ for once-punctured surfaces. In these coordinates the expression for $\Box_i$ is given by
\begin{equation}
\label{equation: box}
\begin{split}
\Box_i=\frac{1}{2} \Big( \dd A_i \wedge \dd B_i + \dd G_i \wedge \dd A_i + \dd B_i \wedge \dd (G_i - A_i) + 2\dd (A_i+B_i) \wedge \dd C_i  \\
+ 2 \dd (G_i - A_i + B_i) \wedge \dd D_i + \dd (G_i - A_i + B_i) \wedge \dd E_i
+ \dd E_i\wedge \dd (A_i + B_i)  \\
+ \dd (A_i+B_i-E_i) \wedge \dd (G_i-A_i+B_i-E_i) + 2\dd (G_i+2B_i-2E_i) \wedge \dd F_i  \Big).
\end{split}
\end{equation}

For the remaining surfaces, the coordinates we use are the weights on the branches $a_i, b_i$ as on \figref{triangle and square} (left) on each of a trigon as on \figref{fig:cactus} and the weights of  branches on the subsurface $C$ chosen as before.

In these coordinates we write the symplectic forms as follows:
\begin{equation}
\label{equation: symp forms rest}
\omega_{0,n} = \omega_{0,5} + \sum_{i=1}^{n-5} \Delta_i, \,\,\,\,\, \omega_{1,n} = \omega_{1,2} + \sum_{i=1}^{n-2} \Delta_i, \,\,\,\,\, \omega_{g,n} = \omega_{g,1} + \sum_{i=1}^{n-1} \Delta_i,  
\end{equation}
where $\Delta_i$ is the contribution of the $3$ switches of the $i-$th trigon as on \figref{fig:cactus}. As before, to write the explicit expression for $\Delta_i$, we use auxiliary coordinates $c_i$, that can be expressed through other variables inductively in the following way:
$$
c_{i+1} = c_i - 2a_i + 2b_i,
$$
with $c_1 = 2( z_1+z_2+z_3+z_4)$ for the punctured spheres, $c_1 = 2(s_1+s_2+s_3-s_4)$ for the punctured tori and $c_1 = y_1-y_2+y_3+y_4$ for the remaining surfaces. In these coordinates the expression for $\Delta_i$ is given by

\begin{equation}
\label{equation: delta}
\Delta_i = 2\dd a_i\wedge\dd b_i + \dd c_i \wedge (a_i-b_i).
\end{equation}
We note that this choice of coordinates allows us naturally embed the space $W(\tau_{g,n})$ into the space $W(\tau_{g,n+1})$ and under this embedding
$$
\omega_{g,n+1}-\omega_{g,n} = \Delta_{k},
$$
with $k=n-4$ for punctured spheres, $k=n-1$ for punctured tori, $k=n$ for the remaining surfaces.

\subsection{Induction by genus.} \label{subsection: ind genus}

In this subsection we run the induction by genus for closed and once-punctured surfaces of genus at least $2$. We use the recurrent expressions for $\omega_{g+1,0}$ and $\omega_{g+1,1}$ from \eqnref{eq: recursive puncture}
and the induction hypothesis to compute their top wedge powers. 

In the case of closed surfaces we obtain
$$
\omega_{g+1,0}^{3g} = (\omega_{g,0}+\Box_{g-1})^{3g} =\sum_{i=0}^{3g} \binom{3g}{i} \omega_{g,0}^i\wedge\Box_{g-1}^{3g-i}.
$$
Notice that $\omega_{g,0}^k=0$ for $k>3g-3$, since $\omega_{g,0}$ depends only on $6g-6$ coordinates, and also $\Box_{g-1}^k=0$ for $k>3$ since the form $\Box_{g-1}$ depends only on $7$ coordinates, see the \eqnref{equation: box}. So the only non-zero monomial in the above expression is $\binom{3g}{3} \omega_{g,0}^{3g-3}\wedge\Box_{g-1}^{3}$, and therefore we get:
$$
\frac{1}{(3g)!}\omega_{g+1,0}^{3g} = \frac{1}{(3g)!} (\omega_{g,0} + \Box_{g-1})^{3g} = \frac{1}{(3g)!} \binom{3g}{3}  \omega_{g,0}^{3g-3} \wedge\Box_{g-1}^3= \frac{1}{(3g-3)!} \omega_{g,0}^{3g-3}\wedge \frac{\Box_{g-1}^3}{6}.
$$
A direct calculation shows that
\begin{equation}
\label{equation: box cube}
\Box_i \wedge \Box_i \wedge \Box_i  = \pm \, 24 \dd A_i \wedge \dd B_i \wedge \dd C_i \wedge \dd D_i \wedge \dd E_i \wedge \dd F_i + \dd G_i \wedge [\dotsc].
\end{equation}
Recall that the coordinate $G_i$ is a linear combination of the coordinates $x_1, x_2, x_3$ and $B_j, E_j$ for $j<i$, so the product $\omega_{g,0}^{3g-3}\wedge \dd G_{g-1}$ vanishes. Having this in mind we obtain:
\begin{equation*}
\frac{1}{(3g)!}\omega_{g+1,0}^{3g}=\frac{1}{(3g-3)!} \omega_{g,0}^{3g-3}\wedge \frac{\Box_{g-1}^3}{6}= \\
=\pm \, \frac{1}{(3g-3)!}\omega_{g,0}^{3g-3}\wedge \frac{24}{6} \dd A_{g-1} \wedge \dd B_{g-1} \wedge \dd C_{g-1} \wedge \dd D_{g-1} \wedge \dd E_{g-1} \wedge \dd F_{g-1}.
\end{equation*}
Since the lattice $\Lambda(\tau_{g+1,0})$ coincides with the standard lattice in the coordinates  $\{ x_i,A_i,B_i,C_i,D_i,E_i,F_i \}$, by the calculation above and the induction hypotheses the ratio between measures $\mu_{\omega}$ and $\mu_{\ZZ}$ equals to
$$\frac{24}{6} 2^{2g-3} = 2^{2(g+1)-3},$$
as required.

For the once-punctured surfaces, we proceed similarly. First we notice that the only non-zero monomial of the sum
$$
\frac{1}{(3g+1)!}\omega_{g+1,1}^{3g+1} = \frac{1}{(3g+1)!} (\omega_{g,1} + \Box_{g-1})^{3g+1} = \frac{1}{(3g+1)!}\sum_{i=0}^{3g+1} \binom{3g+1}{i} \omega_{g,1}^i\wedge\Box_{g-1}^{3g+1-i}
$$
is $\binom{3g+1}{3}  \omega_{g,1}^{3g-2}\wedge \Box_{g-1}^3$, and therefore
$$
\frac{1}{(3g+1)!}\omega_{g+1,1}^{3g+1}= \frac{1}{ (3g-2)!} \omega_{g,1}^{3g-2}\wedge \frac{\Box_{g-1}^3}{6}.
$$
Again, as in the previous case, the product $\omega_{g,1}^{3g-2}\wedge\dd G_i$ vanishes, and therefore we obtain
$$
\frac{1}{(3g+1)!}\omega_{g+1,1}^{3g+1}= \frac{1}{ (3g-2)!} \omega_{g,1}^{3g-2}\wedge \frac{24}{6} \dd A_{g-1} \wedge \dd B_{g-1} \wedge \dd C_{g-1} \wedge \dd D_{g-1} \wedge \dd E_{g-1} \wedge \dd F_{g-1}.
$$
The lattice $\Lambda(\tau_{g+1,1})$ coincides with a sublattice of the standard lattice in the coordinates  $\{ y_i,A_i,B_i,C_i,D_i,E_i,F_i \}$, given by the condition that the sum $y_1-y_2+y_3+y_4$ is even (see \secref{subsection: g21}). Then by the calculation above and the induction hypotheses, the ratio between measures $\mu_{\omega}$ and $\mu_{\ZZ}$ equals to
$$\frac{24}{6} 2^{2g-2} = 2^{2(g+1)-2},$$
as required.

\subsection{Induction by the number of punctures.}
\label{sec: ind punc}
In this subsection we run the induction by the number of punctures for the remaining surfaces using the same approach as in the \secref{subsection: ind genus}.
We use the recurrent expressions for $\omega_{0,n+1}$, $\omega_{1,n+1}$ and $\omega_{g,n+1}$ from \eqnref{equation: symp forms rest}
and the induction hypothesis to compute their top wedge powers. 
In the case of punctured spheres we obtain
$$  \omega_{0,n+1}^{n-2} =  (\omega_{0,n} +\Delta_{n-4})^{n-2} =  \sum_{i=0}^{n-2} \binom{n-2}{i} \omega_{0,n}^i\wedge\Delta_{n-4}^{n-2-i}.$$
Notice that $\omega_{0,n}^k=0$ for $k>n-3$, since $\omega_{0,n}$ depends only on $2n-6$ coordinates, and also $\Delta_{n-4}^k=0$ for $k>1$ since the form $\Delta_{n-4}$ depends only on 3 coordinately, see the \eqnref{equation: delta}. Thus the only non-zero monomial in the above expression is $\binom{n-2}{n-3} \omega_{0,n}^{n-3}\wedge\Delta_{n-4}$, and therefore we get:  
$$ \frac{1}{(n-2)!} \omega_{0,n+1}^{n-2} =  \frac{1}{(n-3)!}   \omega_{0,n}^{n-3}\wedge\Delta_{n-4}.$$
Recall that the coordinate $c_i$ is a linear combination of the coordinates $z_1,\dotsc,z_4$ and  $a_j,b_j$ for $j<i$, so the product $ \omega_{0,n}^{n-3}\wedge \dd c_{n-4}$ vanishes. We obtain
$$\frac{1}{(n-2)!} \omega_{0,n+1}^{n-2}
= \frac{1}{(n-3)!}   \omega_{0,n}^{n-3}\wedge\Delta_{n-4}
=\frac{1}{(n-3)!}   \omega_{0,n}^{n-3}\wedge
2 \dd a_{n-4}\wedge\dd b_{n-4}.$$ 
Since the lattice $\Lambda(\tau_{0,n+1})$  coincides with the standard lattice in the coordinates $\{a_i,b_i,z_i\}$,
by the calculation above and the induction hypotheses the ratio between measures $\mu_{\omega}$ and $\mu_{\ZZ}$ equals to
$$ 2 \cdot 2^{n-3} = 2^{(n+1)-3},$$
as required.

In the case of punctured tori we similarly get
$$ \frac{1}{(n+1)!} \omega_{1,n+1}^{n+1} = \frac{1}{(n+1)!} (\omega_{1,n} +\Delta_{n-1})^{n+1} = \frac{1}{n!} \omega_{1,n}^{n} \wedge \Delta_{n-1}=\frac{1}{n!} \omega_{1,n}^{n} \wedge
2 \dd a_{n-1}\wedge\dd b_{n-1}. $$
Since the lattice $\Lambda(\tau_{1,n+1})$ coincides with the standard lattice in the coordinates $\{a_i,b_i,s_i\}$,
by the calculation above and the induction hypotheses the ratio between measures $\mu_{\omega}$ and $\mu_{\ZZ}$ equals to
$$ 2 \cdot 2^{n-1} = 2^{n},$$
as required.

Finally, for the punctured surfaces of genus at least $2$ we obtain
\begin{equation}
\begin{split}
\frac{1}{(3g-2+n)!} \omega_{g,n+1}^{3g-2+n} = \frac{1}{(3g-2+n)!} (\omega_{g,n} +\Delta_n)^{3g-2+n} =\\ \frac{1}{(3g-3+n)!} \omega_{g,n}^{3g-3+n} \wedge \Delta_n=\frac{1}{(3g-3+n)!} \omega_{g,n}^{3g-3+n} \wedge
2 \dd a_{n}\wedge\dd b_{n}.
\end{split}    
\end{equation}

The the lattice $\Lambda(\tau_{g,n+1})$
coincides with the standard lattice in the coordinates $\{a_i,b_i,  y_i,A_i,B_i,C_i,D_i,E_i,F_i\}$, such that the sum $y_1-y_2+y_3+y_4$ is even (see \secref{subsection: g21}). Therefore by the calculation above and the induction hypotheses the ratio between measures $\mu_{\omega}$ and $\mu_{\ZZ}$ equals to
$$ 2 \cdot 2^{2g+n-3} = 2^{2g+(n+1)-3},$$
as required. 

The computations of the \secref{subsection: ind genus} and \secref{sec: ind punc} together with base cases considered in \secref{subsection: base cases} complete the proof of the main theorem.
\begin{flushright}{$\square$}\end{flushright}


\begin{thebibliography}{[LMir08]}
\bibitem[Ar19]{Arana-Herrera} F. Arana-Herrera, {\em Normalization of Thurston measures on the space of measured geodesic laminations}, In preparation. 


\bibitem[BS\"o01]{BS}
F. Bonahon and Y. S\"ozen, {\em The Weil-Petersson and Thurston Symplectic forms}, Duke Math. J. 108 (2001), no. 3, 581--597.

\bibitem[DGZZ19]{DGZZ19} 
V. Delecroix, E. Goujard, P. Zograf, and A. Zorich {\em Masur-Veech volumes, frequencies of simple closed geodesics and intersection numbers of moduli spaces of curves}, In preparation.


\bibitem[Mas85]{Masur} H. Masur, {\em Ergodic actions of the mapping class group}, Proc. AMS 94, (1985).


\bibitem[Mir08]{Maryam Ergodic} Maryam Mirzakhani, {\em Ergodic Theory of the Earthquake Flow}, International Mathematics Research Notices, Volume 2008

\bibitem[RSo17]{Currents} Kasra Rafi, Juan Souto, {\em Geodesics Currents and Counting Problems}, {\tt arXiv:1709.06834}, (2017). 

\bibitem[PH92]{Penner-Harer}
R. Penner and J. Harer, {\em Combinatorics of Train Tracks}, Princeton University Press, (1992).


\bibitem[Mar16]{Bruno}
Bruno Martelli, {\em An Introduction to Geometric Topology},
{\tt http://people.dm.unipi.it/martelli/Geometric\_topology.pdf}



\end{thebibliography}
\end{document}